\documentclass[12pt]{article}
\usepackage{amsmath,amsthm,amssymb,amsfonts}
\usepackage{latexsym}
\usepackage{psfig}
\usepackage[dvips]{graphicx}

\DeclareMathOperator{\conv}{conv}
\newcommand{\bfa}{{\mathbf a}}
\newcommand{\bfb}{{\mathbf b}}
\newcommand{\bfc}{{\mathbf c}}
\newcommand{\bfp}{{\mathbf p}}
\newcommand{\bfr}{{\mathbf r}}
\newcommand{\bfs}{{\mathbf s}}
\newcommand{\bft}{{\mathbf t}}
\newcommand{\bfu}{{\mathbf u}}
\newcommand{\bfv}{{\mathbf v}}
\newcommand{\bfx}{{\mathbf x}}

\newtheorem{thm}{Theorem}
\newtheorem{lemma}[thm]{Lemma}
\newtheorem{conj}[thm]{Conjecture}

\title{A Degree Bound for \\Codimension Two Lattice Ideals}
\author{Leah H. Gold \\ Texas A \& M University}
\date{October 17, 2002}

\begin{document}
\thispagestyle{empty}

% this defines an arrow that can be made as long as desired

\def \a{{\mathrel{\smash-}}{\mathrel{\mkern-8mu}}
{\mathrel{\smash-}}{\mathrel{\mkern-8mu}}
{\mathrel{\smash-}}{\mathrel{\mkern-8mu}}}

\maketitle

\begin{abstract}
Herzog and Srinivasan have conjectured that for any homogeneous
$k$-algebra, the degree is bounded above by a function of the maximal
degrees of the syzygies. Combining the syzygy quadrangle decomposition
of Peeva and Sturmfels and a delicate case analysis, we prove that
this conjectured bound holds for codimension~2 lattice ideals.
\end{abstract}

\section{Introduction}
\paragraph{} 
Let $R=k[x_1,\ldots,x_n]$ be a polynomial ring in $n$ variables over a
field $k$, let $\deg(x_i) = 1$, and let $I \subset R$ be a homogeneous
ideal. If the Hilbert polynomial of $R/I$ is $\sum_{i=0}^{m} a_i t^i$,
then the {\it degree} of the ideal $I$, written $\deg(I)$, is simply
$a_m \, m!$.

In this section we briefly describe the progress to date on bounding
the degree of an ideal. In particular, we recall several conjectures
which were made about the degree and discuss what is known about the
conjectures. In Section~\ref{latticeSection} we define codimension~2
lattice ideals and explain Peeva and Sturmfels' decomposition of the
resolution of any such ideal.  Finally, in Section~\ref{boundSection}
we use the decomposition and a careful case analysis of the possible
syzygies to prove that the conjectured bound on the degree holds for
codimension~2 lattice ideals.

A resolution is called {\it pure} if at each step there is only a
single degree. That is, the resolution looks like
\[
0  \rightarrow R(-d_p)^{b_p} 
\rightarrow R(-d_{p-1})^{b_{p-1}} 
\rightarrow \dots
\rightarrow R(-d_2)^{b_2} 
\rightarrow R(-d_1)^{b_1} 
\rightarrow R.
\] 

Huneke and M. Miller proved the following formula for the degree
of a Cohen-Macaulay algebra with a pure resolution \cite{HM}.

\begin{thm}[Huneke \& M. Miller] \label{thm1}
  Let $R/I$ be a Cohen-Macaulay algebra with a pure resolution as
  displayed above. Then $ \deg(I) = \frac{\prod_{i=1}^{p} d_i}{p!}$.
\end{thm}

One might hope that when the resolution is not pure that it is
possible to write a similar closed formula for the degree in terms of
the degrees of the syzygies. This does not appear to be the case,
however, Huneke and Srinivasan made a conjecture using similar
formulas to bound the degree \cite{HS}.

\begin{conj}[Huneke \& Srinivasan] \label{conj1}
  Let $R/I$ be a Cohen-Macaulay algebra with resolution of the form
  \begin{center}
  $ 0 \rightarrow \bigoplus\limits_{j \in J_p} R(-d_{p,j}) \rightarrow
  \dots \rightarrow \bigoplus\limits_{j \in J_2} R(-d_{2,j})
  \rightarrow \bigoplus\limits_{j \in J_1} R(-d_{1,j}) \rightarrow R$.
  \end{center}
  Let $m_i = \min \; \{d_{i,j} \in J_i\}$ be the minimum degree shift
  at the $i$th step and let $M_i = \max \;\{d_{i,j} \in J_i\}$ be the
  maximum degree shift at the $i$th step. Then
  \[
  \frac{\prod_{i=1}^p m_i}{p!}  \;\leq \; \deg(I) \;\leq\;
  \frac{\prod_{i=1}^p M_i}{p!}.
  \]
\end{conj}

Notice that since $R/I$ is Cohen-Macaulay, $p$ is the codimension of
$I$.

Due to work by Herzog and Srinivasan \cite{HS}, Conjecture~\ref{conj1}
is known to be true for the following types of ideals.
\begin{list}{ --}{\setlength{\itemsep}{-.2 \baselineskip}}
  \item complete intersections
  \item perfect ideals with quasipure resolutions 
     ($ d_{i,j} \leq d_{i+1,j} \mbox{ for all } i, j$)
  \item perfect ideals of codimension~2
  \item Gorenstein ideals of codimension~3 generated by 5 elements (the
  upper bound holds for all codimension~3 Gorenstein ideals)
  \item perfect stable monomial ideals (as defined by Eliahou \& Kervaire
  \cite{EK})
  \item perfect squarefree strongly stable monomial ideals
  (see Aramova, Herzog, Hibi \cite{AHH})
\end{list}

Generalizing even further, one might want to omit the Cohen-Macaulay
restriction. Consider $I = (x^2, x y) \subset k[x,y]$. Then $\deg(I)=
1$, $m_1 = 2$ and $m_2 = 3$, but $\frac{(2)(3)}{2!} \geq 1$. So we
know that the lower bound does not hold for non-Cohen-Macaulay
algebras and therefore we consider just the upper bound in the
non-Cohen-Macaulay case.

\begin{conj}[Herzog \& Srinivasan]\label{conj2}
  Let $I$ be a homogeneous ideal of codimension~$d$ and $M_i$ as
  defined above, then $\deg(I) \leq \frac{\prod_{i=1}^d M_i}{d!}$.
\end{conj}

Herzog and Srinivasan showed that Conjecture~\ref{conj2} is true 
in three cases.
\begin{list}{ --}{\setlength{\itemsep}{-.2 \baselineskip}
  }
  \item stable monomial ideals (as defined by Eliahou \& Kervaire \cite{EK})
  \item squarefree strongly stable monomial ideals (see Aramova,
   Herzog, Hibi \cite{AHH})
  \item ideals with a $q$-linear resolution (all the generators in
  degree q and all the syzygies are linear)
\end{list}

Prior to the result presented here, the above cases formed a complete
list of all known cases where the conjectures are true.

\section{Codimension~2 lattice ideals}\label{latticeSection}
\paragraph{} 
Lattice ideals are a slight generalization of toric
ideals. Codimension~2 lattice ideals were studied by Peeva and
Sturmfels in their paper [PS]. We briefly describe here the
relevant results from their paper, namely, the construction of an
explicit resolution of any such ideal.

We begin by defining a lattice ideal. Let $R = k[x_1, \ldots, x_n]$ be
a polynomial ring and for any nonnegative integer vector $\bfa = (a_1,
\ldots, a_n)$, let $x^\bfa = x_1^{a_1} \cdots x_n^{a_n}$. For any
lattice ${\cal L} \subset {\mathbb Z}^n$ we define 
\[
 I_{\cal L} \; = \; ( \bfx^{\bfa_{+}} - \bfx^{\bfa_{-}} \; | \; \bfa \in
 {\cal L} )
\] 
\noindent 
where $\bfa_+$ is the positive part of the vector $\bfa$ and $\bfa_-$
is the negative part of $\bfa$.  That is, in the $i$th component,
$(\bfa_+)_i = \bfa_i$ if $\bfa_i \geq 0$ and zero otherwise. We define
${\bfa_-}$ in a similar manner. We consider only lattices with no
nonnegative vectors in order to ensure the lattice ideal is
homogeneous with respect to some positive grading.

We may define a multigrading on $R$, and also on $I_{\cal L}$, by the
group ${\mathbb Z}^n / \cal L$. We will move back and forth between
this grading and the standard grading. It should be clear from context
which one is meant.

The codimension of $I_{\cal L}$ is the minimal number of generators of
the lattice $\cal L$. When $I_{\cal L}$ has codimension~2, Peeva and
Sturmfels constructed a resolution for $I_{\cal L}$ in the following
way.

Let $\bfc$ be a multidegree and let ${\bfx^\bfa}$ be a monomial of
degree $\bfc$. Then there is a correspondence between monomials of
degree $\bfc$ and vectors $\bfu \in {\mathbb Z}^2$ such that ${\mathrm
B}\bfu \leq \bfa$. The monomial ${\bfx^\bfa}-{\mathrm B}\bfu$
corresponds to the vector $\bfu$.  Define the polytope $P_\bfa =
\conv(\{\bfu \in {\mathbb Z}^2 \, | \, {\mathrm B}\bfu \leq \bfa \})$.
Notice that $P_\bfa$ and $P_\bfb$ are lattice translates of each other
if and only if $\bfa - \bfb \in \cal L$. So, we generally write
$P_\bfc$ instead of $P_\bfa$.

Peeva and Sturmfels showed that each multidegree in which there is a
minimal syzygy corresponds to a primitive polytope. In particular,
first syzygies correspond to line segments, second syzygies correspond
to triangles and third syzygies correspond to quadrangles. Further,
the syzygy triangles consist of three syzygy line segments and the
syzygy quadrangles consist of four syzygy triangles. For details on
this correspondence, see Peeva and Sturmfels paper \cite{PS}. A
resolution of the ideal generated by the binomials corresponding to
the four segments is found by the following method.

Let $P_\bfc$ be a polytope corresponding to a syzygy quadrangle. We
start by writing the two generators of $I_{\cal L}$ corresponding to
the sides of the quadrangle as $\alpha = \alpha' - \alpha''$ and
$\beta = \beta' - \beta''$. Then we determine vectors $\mathbf p$,
$\mathbf r$, $\mathbf s$, and $\mathbf t$ by taking the greatest
common divisors of a term of $\alpha$ and a term of $\beta$. For
example, choose $\mathbf p$ such that $x^{\mathbf p} = \gcd(\alpha',
\beta')$. We set the remaining factors to be ${\mathbf x}^{{\bfu_+}}$,
${\mathbf x}^{{\bfu_-}}$, ${\mathbf x}^{{\bfv_+}}$, and ${\mathbf
x}^{{\bfv_-}}$ and so we have $\alpha = {\mathbf x}^{{\mathbf
u_+}}{\mathbf x}^{\mathbf p}{\mathbf x}^{\mathbf t} - {\mathbf
x}^{{\bfu_-}}{\mathbf x}^{\mathbf r}{\mathbf x}^{\mathbf s}$ and
$\beta = {\mathbf x}^{\bfv_+}{\mathbf x}^{\mathbf p}{\mathbf
x}^{\mathbf s} - {\mathbf x}^{\bfv_-}{\mathbf x}^{\mathbf r}{\mathbf
x}^{\mathbf t}$.

A diagonal vector of the quadrangle is a sum or difference of the two
edge vectors. Hence we can derive representations for the generators
which correspond to the diagonals from the generators for $\alpha$ and
$\beta$ by taking the sum or difference of the exponent vectors of the
binomials $\alpha$ and $\beta$.  This procedure gives $\gamma =
{\mathbf x}^{{\bfu_+}} {\mathbf x}^{{\bfv_+}} {\mathbf x}^{2 \mathbf
p} - {\mathbf x}^{{\bfu_-}} {\mathbf x}^{{\bfv_-}} {\mathbf
x}^{2\mathbf r}$ and $\delta = {\mathbf x}^{{\bfu_+}} {\mathbf
x}^{{\bfv_-}} {\mathbf x}^{2 \mathbf t} - {\mathbf x}^{{\bfu_-}}
{\mathbf x}^{{\bfv_+}} {\mathbf x}^{2 \mathbf s}$. Notice that in order
to generate an ideal of codimension~2, the four generators of the
ideal $I_{\cal L}$ cannot share a common factor.

Putting all of this into a sequence, the resolution of the four
generators derived from $P_\bfc$ has the form
\[
0  \rightarrow
 R  
{ {\mathop{\a {\mathrel{\smash-}}
{\mathrel{\mkern-8mu}}\rightarrow}\limits^{
\begin{pmatrix} 
-{\bf x}^{\bf s}\cr {\bf x}^{\bf t}\cr
{\bf x}^{\bf r}\cr  -{\bf x}^{\bf p} 
\end{pmatrix}
} }}
 R^4 
 {\mathop{{\mathrel{\smash-}}{\mathrel{\mkern-8mu}}
\a\a\a\a\a\a\a\a\a\a\a\a\a\a
{\mathrel{\smash-}}{\mathrel{\mkern-8mu}}
\rightarrow}\limits^{
 \begin{pmatrix}
{\bf x}^{{\bf v}_+}{\bf x}^{\bf p}&
{\bf x}^{{\bf v}_-}{\bf x}^{\bf r}
&-{\bf x}^{{\bf v}_-}{\bf x}^{\bf t}&
-{\bf x}^{{\bf v}_+}{\bf x}^{\bf s}\cr
{\bf x}^{{\bf u}_-}{\bf x}^{\bf r}&
{\bf x}^{{\bf u}_+}{\bf x}^{\bf p}
&{\bf x}^{{\bf u}_-}{\bf x}^{\bf s}&
{\bf x}^{{\bf u}_+}{\bf x}^{\bf t}\cr
 -{\bf x}^{\bf t}& -{\bf x}^{\bf s}&
0&0\cr
0&0& {\bf x}^{\bf p}& {\bf x}^{\bf r}   
\end{pmatrix}
 } }
 \, R^4 \,
{\mathop{\a\a{\mathrel{\smash-}}{\mathrel{\mkern-8mu}}
{\mathrel{\smash-}}{\mathrel{\mkern-8mu}}
\rightarrow}\limits^{
{(\alpha\ \  \beta\  \ \gamma \  \ \delta)}}}
R
\]

where 
\[
\begin{array}{ll}
\alpha = \bfx^{\bfu_+}\bfx^\bft{\bfx}^\bfp -
\bfx^{\bfu_-}\bfx^\bfs\bfx^\bfr \;\;\;\;\;& \beta =
\bfx^{\bfv_+}\bfx^\bfs{\bfx}^\bfp -
\bfx^{{\bfv_-}}\bfx^\bft\bfx^\bfr \;\;\;\;\; \\ \gamma =
\bfx^{{\bfu_+}}\bfx^{{\bfv_+}}\bfx^{2 \bfp} -
\bfx^{\bfu_-}{\bfx}^{\bfv_-}\bfx^{2\bfr} \;\;\;\;\;& \delta =
\bfx^{\bfu_+}\bfx^{\bfv_-}\bfx^{2 \bft} -
\bfx^{\bfu_-}\bfx^{\bfv_+}\bfx^{2 \bfs}\;\;\;\;\;
\end{array}.
\]

The resolutions corresponding to the syzygy quadrangles may then be
used to build a resolution for $R/I_{\cal L}$.

\begin{thm}[Peeva \& Sturmfels]\label{PSthm}
  If $R/I_{\cal L}$ is not Cohen-Macaulay, then the sum of the
  complexes corresponding to syzygy quadrangles is a minimal free
  resolution of $R/I_{\cal L}$.
\end{thm}

\section{Bounding the degree for codimension~2 lattice ideals}\label{boundSection}
\paragraph{} 
Using the decomposition of the resolution in terms of the syzygy
quadrangles given in the previous section, we can now show that
Conjecture~\ref{conj2} is true for codimension~2 lattice ideals.

\begin{thm}\label{mythm}
  If $I_{\cal L}$ is a homogeneous codimension~2 lattice ideal, then
 Conjecture~\ref{conj2} holds. That is, for a homogeneous
 codimension~2 lattice ideal $I_{\cal L}$, if $M_1$ is the maximal
 degree of the generators of $I_{\cal L}$ and $M_2$ is the maximal
 degree of the syzygies on the generators, then
  \[
  \deg(I_{\cal L}) \leq \frac{M_1M_2}{2}.
  \]
\end{thm}

In order to prove this theorem, we first prove a special case.

\begin{lemma}\label{mylemma}
  Let $I_{\cal L}$ be a lattice ideal and let $J$ be an ideal whose
  four generators are associated to a single syzygy quadrangle of
  $I_{\cal L}$. Then Theorem~\ref{mythm} holds for $J$.
\end{lemma}

\begin{proof}
Let $J$ be the ideal whose four generators are associated to a single
syzygy quadrangle of $I_{\cal L}$. Using the resolution described in
the previous section, we can write down the Hilbert series for
$R/J$. That is, 
\[
H_{R/J}(y) = \frac{f(y)}{(1-y)^n}
\, \mbox{ where } \,
f(y)= \sum_{i=1}^n \sum_{j \in J_i} (-1)^i y^{d_{i,j}}. 
\]

Canceling powers of $(1-y)$, we obtain 
\[
H_{R/J}(y) = \frac{g(y)}{(1-y)^{n-2}}.
\]
So $\deg(J) = g(1)= \frac{1}{2} f''(1)$.

For any vector $\bfv=|v|=(v_1,v_2,\ldots, v_n)$, let $v = v_1 +
v_2 + \cdots + v_n$. Using this notation and our knowledge of the
$d_{i,j}$ from the resolution, we can write $\deg(J)$ in terms of
$u_+$, $u_-$, $v_+$, $v_-$, $p$, $r$, $s$, and $t$.

Since $\alpha$ and $\beta$ are homogeneous polynomials, there are
relations between the eight variables $u_+$, $u_-$, $v_+$, $v_-$, $p$,
$r$, $s$, and $t$ which arise because the degrees of the terms in the
binomials are equal.  Using these relations, we can eliminate $u_{-}$
and $v_{-}$ and write $\deg(J)$ in terms of the other six variables,
$u_+$, $v_+$, $p$, $r$, $s$, and $t$. So, 
\[ 
\deg(J) = u_{+} v_{+} + u_{+} p + v_{+} p + p^2 - p r + u_{+}s + p s
+v_{+}t +p t. \]

Now, what are the possibilities for $M_1$ and $M_2$? $M_1$ could be
$\deg(\alpha)$, $\deg(\beta)$, $\deg(\gamma)$, or $\deg(\delta)$ and
$M_2$ could be $\deg(\gamma)+s$, $\deg(\gamma)+t$, $\deg(\delta)+p$,
or $\deg(\delta)+r$. We proceed by investigating these cases.

We begin by using the fact that the syzygies are homogeneous to
describe some relations on the exponents.
\[
\begin{array}{c}
v_{+} + p + \deg(\alpha) = u_{-} + r + \deg(\beta)  = t + \deg(\gamma)\\
v_{-} + r + \deg(\alpha) = u_{+} + p + \deg(\beta) = t + \deg(\gamma)\\
v_{-} + t + \deg(\alpha) = u_{-} + s + \deg(\beta) = p + \deg\delta)\\
v_{+} + s + \deg(\alpha) = u_{+} + t + \deg(\beta) = r + \deg(\delta)\\
\end{array}
\]

From these equalities, we can distill the inequalities $\deg(\gamma) +
\deg(\delta) \geq 2 \deg(\alpha)$ and $\deg(\gamma) + \deg(\delta)
\geq 2 \deg(\beta)$. Hence $M_1 = \deg(\gamma) \mbox{ or }
\deg(\delta)$. Since $\gamma$ and $\delta$ are interchangeable, we can
assume $M_1 = \deg(\delta)$.

This leaves us with four cases to check corresponding to the four
possible values of $M_2$. In each case we consider the expression for
$M_1M_2 - 2 \deg(J)$. We expand the expression in terms of $u_+$,
$u_-$, $v_+$, $v_-$, $p$, $r$, $s$ and $t$ and then eliminate two of
the variables using the equations arising from the homogeneity
conditions. The choice of which variables to eliminate is not obvious,
but there is always a nice choice which makes it easier to show that
the expression is nonnegative. Then, in each case, we can show that
the expression is nonnegative by using the inequalities that arise
from the choices of $M_1$ and $M_2$. Which inequalities were necessary
and how to use them were not obvious at first glance so a computer
program PORTA \cite{P} was used to help reduce the inequalities. Once
it was clear what we should look for, it was easy to do these by hand.

Consider the case where $M_2 = \deg(\gamma) + t$. If we eliminate
$u_{-}$ and $v_{+}$, the expression for $M_1 M_2 - \deg(J)$ can be
rewritten as

$u_{+}^2 + v_{-}^2 + u_{+}(p-s) + v_{-}(r-s) + u_{+}(t-r) + u_{+}t +
v_{-}(t-p) + v_{-}t + 2t^2$.

The choices of $M_1$ and $M_2$ in this case imply that $p \geq s$, $r
\geq s$, $t \geq r$, and $t \geq p$. It is clear, therefore, that the
above expression is nonnegative.

The other three cases look similar although different eliminations and
inequalities are used for each case. 

So, for every possible choice of $M_1$ and $M_2$, the expression
$M_1M_2 - 2 \deg(J)$ is nonnegative.  Therefore all ideals $J$ arising
from a syzygy quadrangle satisfy the bound of Conjecture~\ref{conj2}.
\end{proof}

\begin{proof}[Proof of Theorem~\ref{mythm}.]
If $R/I_{\cal L}$ is Cohen-Macaulay, then we know from Herzog and
Srinivasan's paper \cite{HS} that it satisfies the bound. So, let us
assume $R/I_{\cal L}$ is not Cohen-Macaulay.  We may construct a
resolution for $R/I_{\cal L}$ via its syzygy quadrangles.

Let $J$ be an ideal whose four generators are associated to a single
syzygy quadrangle of $I_{\cal L}$. Since, according to
Theorem~\ref{PSthm}, the syzygies from this resolution are also
syzygies of $R/{\cal L}$, we know that $M_1(J) \leq M_1(I_{\cal L})$
and $M_2(J) \leq M_2(I_{\cal L})$. Together these imply that
\[ 
\frac{M_1(J)M_2(J)}{2} \leq \frac{M_1(I_{\cal L})M_2(I_{\cal L})}{2}.
\]
On the other hand, since $J \subset I_{\cal L}$, we have that
$\deg(I_{\cal L}) \leq \deg(J)$.

By Lemma~\ref{mylemma}, we know that the bound holds for $J$. That is
\[
\deg(J) \leq \frac{M_1(J)M_2(J)}{2}.
\]
 Hence,
\[
\deg(I_{\cal L}) \leq \deg(J) \leq \frac{M_1(J)M_2(J)}{2}
\leq \frac{M_1(I_{\cal L})M_2(I_{\cal L})}{2}. \;\;\; 
\]
\end{proof}

\section{Further thoughts}
\paragraph{} 
We have shown that the conjecture of Herzog and Srinivasan is true for
codimension two lattice ideals. For non-Cohen-Macaulay codimension~2
lattice ideals, this bound cannot be tight, that is, we cannot force the
expression $M_1 M_2 - 2 \deg(J)$ to be zero. Suppose we try to force
the expression to be zero, then the squared variables would have to be
zero. In the case where $M_2 = \deg(\gamma) + t$ for instance, it would
force $u_{+}$, $v_{-}$, and $t$ to be zero. So $M_2 = \deg{\gamma}$
which forces $p=r=s=t=0$. Hence $I_{\cal L}$ is Cohen-Macaulay and we
have a contradiction. The other cases are similar.

Although we cannot find an ideal where equality holds, one way to try
to make the above expression small, that is, to make $M_1 M_2$ close
to $2 \deg(J)$, is to choose an ideal where $p=r=s=t$. Doing this also
forces $u_+=u_-$ and $v_+=v_-$. Hence the expression $M_1 M_2 - 2
\deg(J) = {u_+}^2 + {v_+}^2 +u_+p+u_-p+2p^2$ for all choices of $M_1$
and $M_2$. If we then let $u_+=v_+=0$ and $p=1$, we get some ideal of
degree 2 with $M_1=2$ and $M_2=3$. For example, in four variables the
lattice generated by $(1,-1,-1,1)$ and $(1,-1,1,-1)$ gives the ideal
$(a d-b c, a c - b d, a^2 - b^2, c^2 - d^2)$ which has this form. Thus
the bound is quite close to being tight. On the other hand, if we do
not require $u_+ = v_+ = 0$, the expression $M_1 M_2 - 2 \deg(J)$
increases like $(\deg(\alpha))^2$ as $u_+$, $u_-$, or $p$ increases.

The general form of a resolution for codimension~3 or higher lattice
ideals is unknown. These higher codimension lattice ideals do not seem
to lend themselves to the same sort of decomposition as in
codimension~2, so extending the method we used here to prove the bound
for codimension~3 does not seem promising. The case of
non-Cohen-Macaulay ideals of codimension~2 other than lattice ideals
is also still open.

The author thanks Mike Stillman, Irena Peeva and Hal Schenck for all
their help and encouragement.

\footnotesize

%  ****************************************
%  *         BIBLIOGRAPHY                 *
%  ****************************************
%
% These have been checked with the original published source.

\end{document}